\documentclass[12pp]{article}
\usepackage{latexsym,amssymb,amsmath,amscd,amscd,amsthm,amsxtra}
\usepackage[dvips]{graphicx}

\newtheorem{thm}{Theorem}[section]

\newtheorem{cor}[thm]{Corollary}
\newtheorem{pro}[thm]{Proposition}
\newtheorem{ex}[thm]{Example}
\newtheorem{rmk}[thm]{Remark}

\newcommand{\px}{\frac{\partial }{\partial x}}
\newcommand{\py}{\frac{\partial }{\partial y}}
\newcommand{\pz}{\frac{\partial }{\partial z}}

\newcommand{\lon }{\,\rightarrow\,}

\newcommand{\pf}{\noindent{\bf Proof.}\ }

\newcommand{\frkg}{\mathfrak g}

\def\gpd{\,\lower1pt\hbox{$\longrightarrow$}\hskip-.24in\raise2pt
         \hbox{$\longrightarrow$}\,}

\def\qed{\hfill ~\vrule height6pt width6pt depth0pt}

\newcommand{\ZZ}[2]{\rule[#1]{0pt}{#2}}

\begin{document}

\title{{ Linear Poisson Structures on $
\mathbb R^4$
\thanks{ Research partially supported by NSF of China  and the Research Project of ``Nonlinear Science".}}}
%\thanks{1991 \textbf{ Mathematics
%Subject Classification.} Primary 58F05. Secondary 22E65, 22A22,
%58B25,
%58H05.}
\author{ Yunhe SHENG  \\
   Department of Mathematics and LMAM\\ Peking University,
Beijing 100871, China\\
          {\sf email: syh@math.pku.edu.cn} }

\date{}

%\footnotetext{{\it{MSC}}: Primary 58F05. Secondary 17B66, 22E65.}

\maketitle
\begin{abstract}
 We classify all of the
         4-dimensional linear Poisson structures of which the corresponding Lie algebras can be considered as
the extension by a derivation of  3-dimensional unimodular Lie
         algebras. The affine
         Poisson structures on $\mathbb{R}^3$ are totally classified.
\end{abstract}

%\tableofcontents

\section{Introduction}

Linear Poisson structures are in one-to-one correspondence with Lie
algebra structures and usually called Lie-Poisson structures. It is
the most basic and important Poisson structure both for its
exquisite algebraic and geometric properties and its far and wide
applications in physics and other fields of mathematics. In
\cite{LLS}, the authors have classified linear Poisson structures on
$\mathbb R^3$, i.e., give the classification of Lie algebras on
$\mathbb R^3$, see also \cite{Jac} for more details which gives the
classification of 3-dimensional Lie algebras on algebraic closed
field. The idea of using linear Poisson structures to understand the
structures of  Lie algebras can be traced back to the work of Lie.
In this spirit, there have been some suggestions of pursuing this
geometric approach for Lie algebra structures(e.g., see
\cite{Marmo2} and \cite{Marmo1}).

A natural  problem is to classify linear Poisson structures on
$\mathbb R^4$. We find that any 4-dimensional Lie algebra is the
extension of some unimodular 3-dimensional Lie algebra by the
viewpoint of Poisson geometry, so based on the results of the
classification of 3-dimensional Lie algebras in \cite{LLS} and after
the computation of cohomology groups, linear Poisson structures on
$\mathbb R^4$ are totally classified. Furthermore, as mentioned in
\cite{V}, the affine Poisson structures are in one-to-one
correspondence with central extension of the corresponding Lie
algebras and the affine Poisson structures on $\mathbb R^3$ are
totally classified by the way. At last, we give an example of Jacobi
manifold of which the leaves enjoy conformal symplectic structure.

The paper is organized  as follows: In Section 2 we briefly review
the decomposition of linear Poisson structures on $\mathbb R^n$ that
was done in \cite{LLS} and concentrate on $\mathbb R^3$ and $\mathbb
R^4$, obtain the result that any 4-dimensional Lie algebra is the
extension of some unimodular 3-dimensional Lie algebra. In Section 3
we list some useful results that related with the classification of
linear Poisson structures on $\mathbb R^3$ which will be used when
we consider the classification of Lie-Poisson structures on $\mathbb
R^4$. In Section 4 we give a detail description of cohomology groups
with coefficients in trivial representation and adjoint
representation of 3-dimensional Lie algebras. In the last section we
firstly consider all possible extensions of 3-dimensional Lie
algebras and then linear Poisson structures on $\mathbb R^4$ are
totally classified. By the way, we obtain the result of the
classification of affine Poisson structures $\mathbb R^3$.

{\bf Acknowledgement:} ~~ The author would like to give warmest
thanks to Prof. Zhang-ju Liu for his advice and also give thinks to
Prof. G. Marmo for useful comments.

\section{The Decomposition of Lie-Poisson Structures}

In this section, we firstly review the decomposition of linear
Poisson structures on $\mathbb R^n$ \cite{LLS} and based on this, we
mainly concentrate on Lie-Poisson structures on $\mathbb R^3$ and
$\mathbb R^4$. Consider Lie algebras corresponding to them, we find
that every 4-dimensional Lie algebra is the extension, central
extension or extension by a derivation, of some unimodular
3-dimensional Lie algebra.

Throughout the whole paper, $\frak g$ will be a $n$-dimensional Lie
algebra and $\frak g^*$ its dual space with Lie-Poisson structure
$\pi_{\frak g}$ on it. $e_1,\cdots,e_n$ are the basis of Lie algebra
$\frak g$, the corresponding coordinate functions are
$x_1,\cdots,x_n$ and $e^1,\cdots,e^n$ are the dual basis of $\frak
g^*$, the corresponding coordinate functions are $x^1,\cdots,x^n$.
In \cite{LLS}, the authors have given the decomposition of
Lie-Poisson structures on $\mathbb R^n$, Let's recall them briefly.

Let  $\Omega = dx_{1}\wedge dx_{2}\cdots\wedge dx_{n}$ be the
canonical volume
        form  on $\mathbb{R}^n$. Then  $\Omega$ induces an
isomorphism  $\Phi $ from the space of all $i$-multiple vector
fields to the space of all $(n-i)$-forms. Let $d$ denote the usual
exterior differential on forms and
 $$D=(-1)^{k+1}\Phi^{-1} \circ d \circ  \Phi: ~~\mathcal{X}^{k}(\mathbb{R}^n) \lon
\mathcal{X}^{k-1}(\mathbb{R}^n),$$ its pull back under the
isomorphism $\Phi$,
 where $\mathcal{X}^{k}(\mathbb{R}^n)$ denotes the space of all $k$-multiple
vector fields on $\mathbb{R}^n$. An important property of $D$ is
that the Schouten bracket can be written in terms of this operator
as follows \cite{Ko}:
\begin{equation}
\label{eq:schouten} [U,\ V]=D(U\wedge V)-D(U)\wedge V-(-1)^{i}
U\wedge D(V),
\end{equation}
for all $U\in \mathcal{X}^{i}(\mathbb{R}^n)$ and $V \in
\mathcal{X}^{j}(\mathbb{R}^n)$.
       It is obvious that  there is a one-to-one
correspondence between matrices in $\mathfrak{gl(n)}$ and linear
vector fields on $\mathbb{R}^n$, i.e.,
\begin{equation}\label{eq Ahat} A =(a_{ij}) \longleftrightarrow \hat{A}=\sum_{ij}a_{ij}x_{j}\frac{\partial}{\partial
x_{i}} , ~~~~~~ \, ~~~ div_{\Omega}\hat{A}= D(\hat{A}) = \mbox{tr}
A.
\end{equation}
 Moreover, a vector $k \in \mathbb{R}^n$ corresponds to a
constant vector field $\hat{k}$ by translation on $\mathbb{R}^n$ and
satisfies
\begin{equation}\label{eq Ak}
div_{\Omega}\hat{k}= D(\hat{k}) = 0, \, ~~~~~~~~~\,
~~~~~~[\hat{A},~\hat{k} ] = -\hat{Ak},~~~~~~~~~\,~~~~~~\forall ~A
\in \mathfrak{gl(n)}.
\end{equation}

 For  a given  Poisson tensor $\pi$, let
$D(\pi )$ be its modular vector field (see \cite{We}, which is also
called the curl vector field in \cite{DH}). Such a vector field is
always compatible with $\pi$, i.e.,
\begin{equation}
L_{D(\pi )}\pi= [D(\pi ),\ \pi]=0.
\end{equation}
 A Poisson  structure is called unimodular if
$ D(\pi ) = 0$. For a linear Poisson structure $\pi$ on $\mathbb
R^n$, there exists some $k\in \mathbb R^n$ such that $ D(\pi ) =
\widehat{k}$ and  $k$ is also called the modular vector of $\pi$. In
fact $k$ is        always invariant if one takes a
        different volume form. As mentioned in \cite{We}, $k$
        is the modular character of Lie algebra $\mathfrak g$ which corresponds to the linear Poisson structure
        defined as a vector in $\mathfrak g^*$ such that
        $$ \langle k,\xi \rangle =   tr\circ ad(\xi),~~~~~~~~~~~ \,~~~~~~~~~~~~~~~~~
        \forall ~\xi  \in \mathfrak g.$$

        \begin{thm}\label{Th. decom}{\rm\cite{LLS}}
            Any Lie-Poisson structure $\pi_{\frkg}$ on $ {\mathfrak g}^* \cong \mathbb{R}^n $ has a unique
            decomposition:\\
            \begin{equation}\label{linear}
            \pi_{\frkg} = \frac{1}{n-1}\hat{I}\wedge\hat{k} + \Lambda_{{\mathfrak g}},
            \end{equation}
            where $k \in  \mathbb R^n$ is
            the modular vector of $\pi_{\frak g}$ and $\Lambda_{{\frak g}}$ is a linear bi-vector field satisfying
 $D(\Lambda_{{\frak g}})=0$ and                        $ ( \frac{1}{n-1}\hat{k}, \Lambda_{{\frak g}})$
                            is a Jacobi  structure.

 Conversely,  any  such a pair satisfying the  above compatibility conditions defines a  Lie-Poisson structure by
            Formula (\ref{linear}).
        \end{thm}
\begin{rmk}
For more information about Jacobi structures, please see
\cite{ki:local} and  \cite{li:varietes}.
\end{rmk}

Next we give some language of cohomology groups and extension of a
Lie algebra which will be often used later and then based on Theorem
\ref{Th. decom}, we obtain the first result
 that every 4-dimensional Lie algebra is the extension  of some unimodular
3-dimensional Lie algebra.

Recall that for any Lie algebra $\frak g$ and its representation
$\rho:\frak g\longrightarrow \frak gl(V)$ on a vector space $V$, we
have the standard {\em Chevalley cochain complex}
$C^k=Hom(\wedge^k\frak g,V)$ and the coboundary operator
$\delta^k:C^k\longrightarrow C^{k+1}$ is given by
\begin{eqnarray}
\nonumber(\delta
f)(\xi_1,\cdots,\xi_{k+1})&=&\sum^{k+1}_{i=1}(-1)^i(\rho\xi)f(\xi_1,\cdots,\widehat{\xi_i},\cdots,\xi_{k+1})\\\nonumber
&+&\sum_{i<j}(-1)^{i+j}f([\xi_i,\xi_j],\xi_1,\cdots,\widehat{\xi_i},\cdots,\widehat{\xi_j},\cdots,\xi_{k+1}),~
\forall~ f\in C^k,~\xi_1,\cdots,\xi_{k+1}\in \frak g.
\end{eqnarray}
 In particular,
$\delta^0:C^0=V\longrightarrow C^1 $ is given by
\begin{equation}
\delta^0 v(\xi)=(\rho\xi)(v),\quad \forall~ v\in V,~\xi\in \frak g.
\end{equation}
There are two natural representations of $\frak g$ which are trivial
representation on $\mathbb R$ and adjoint representation on itself.

The expression ``$\frak h$ is the {\em central extension} of $\frak
g$ by $\mathbb R$" means that one has a well defined exact sequence
of Lie algebras
\[\begin{CD}
0\longrightarrow\mathbb R@>\iota>> \frak h@>\kappa>>\frak
g\longrightarrow0,
\end{CD}\]
where $\iota(\mathbb R)$ belongs to the center of $\frak h$. This
shows that $\frak h=\frak g\oplus\mathbb R$ as a vector space, and
we must have
$$
[\xi\oplus t,\eta\oplus s]_{\frak
h}=[\xi,\eta]\oplus\omega(\xi,\eta),
$$
where $\omega$ is a {\em Chevalley-Eilenberg} 2-cocycle of $\frak
g$. Furthermore, there is also a one-to-one correspondence between
{\em affine Poisson structures(or modified Lie-Poisson structures)}
on $\frak g^*$ and central extensions of $\frak g$. See more details
in \cite{V}. For convenient, denote by $\frak g_\omega$ the
extension of $\frak g$ decided by $\omega$.

 %Given a unimodular 3-dimensional Lie algebra $\frak g$ with bracket
 %$[\cdot,\cdot]$ and a 2-cocycle $\omega : \frak g\times \frak
 %g\rightarrow \mathbb{R}$, we can define its {\em central extension},
 %denoted by $\frak g_\omega$, $\frak g_\omega=\frak g\oplus
 %\mathbb{R}c$, which is a 4-dimensional Lie algebra with bracket
 %$[\cdot,\cdot]_\omega=[\cdot,\cdot]+\omega(\cdot,\cdot)c$ and $c$ is
 %the new Lie algebra's center.

 Given a  3-dimensional Lie algebra $\frak g$ with bracket
 $[\cdot,\cdot]$ and a derivation $D$, we can define a new
 4-dimensional Lie algebra as the extension of $\frak g$ by the
 derivation $D$({\em D-extension} for convenient), denoted by $\frak
 g_D$, $\frak g_D=\frak g\oplus \mathbb{R}e$, with bracket
 $[\cdot,e]_D=D(\cdot)$.

It is known that if $\omega_1-\omega_2$ is exact, $\frak
g_{\omega_1}$ is
  isomorphic to $\frak g_{\omega_2}$; If $D_1-D_2$ is exact, $\frak g_{D_1}$
    is isomorphic to $\frak g_{D_2}$.

\begin{pro}\label{extension}%%%%%%%%%%%%%%%%%%%%%%%%%%定理1.1
 ~Any {\rm4}-dimensional Lie algebra is the extension, central extension
 or D-extension, of some unimodular {\rm3}-dimensional Lie algebra.
 \end{pro}

  \pf~As stated in the Theorem \ref{Th. decom},
  the Lie-Poisson structure on $\frak g^*$ has the decomposition
 $\pi_{\frak g}=\frac{1}{3}\widehat{I}\wedge\widehat{k}+\Lambda_{\frak g}$. First we
 consider the case that $k=0$.

  $\pi_{\frak g}=\Lambda_{\frak g}$ is the Lie-Poisson
 structure on $\frak g^\ast$, so $[\Lambda_{\frak g}, \Lambda_{\frak g}]=0$.
 Using Formula (\ref{eq:schouten}), we can easily get this equals
 $D(\Lambda_{\frak g}\wedge \Lambda_{\frak g})=0$ and this equals
 $\Lambda_{\frak g}\wedge
 \Lambda_{\frak g}=0$.

 Assume
 $$
 \Lambda_{\frak g}=\sum_{i,j,k=1}^4C^k_{ij}x^k\frac{\partial}{\partial x^i}\wedge
 \frac{\partial}{\partial x^j}=\sum^4_{l=1}x^l\pi_l,
 $$ where
 $\pi_l=\sum_{i,j}C^l_{ij}\frac{\partial}{\partial
 x^i}\wedge\frac{\partial}{\partial x^j}$.  We have $$\Lambda_{\frak
 g}\wedge \Lambda_{\frak g}=0\Longleftrightarrow x^l\pi_l\wedge
 x^l\pi_l=0\Longleftrightarrow \pi_l\wedge \pi_h=0,$$ for any
 $l,h=1,2,3,4$. So $\pi_l$ decide a 2-dimensional subspace $P_l$ and
 $P_l$, $P_h$ intersect a line. If there exists some $\pi_l=0$ or some
 $P_l$ is the linear composition of the others, the problem is easier
 and we can omit them.

 Then consider the lines which are the intersection of some two
 subspaces $P_l$ and $P_h$, there may be three cases below:

 (1).~There are two different subspaces at least and all the subspace
 intersect only one line $L$;

 (2).~There are three linear independent lines $L_1,~L_2,~L_3;$

 (3).~All $P_l$ coincide.

 As in the Case (1), let $L=\mathbb{R}e^4=\cap P_l\in \frak g^*$ and
 $\pi=X\wedge \frac{\partial}{\partial x^4} $ for some linear vector field $X$ on $\frak g^*$.
Denote $L^0=ker(e^4)\subset \frak g$ and so $\frak g=L^0\oplus
\mathbb{R}e_4$. We will show
 that $L^0$ is a 3-dimensional Abelian idea, and then follows that Lie
 algebra $\frak g$ is the D-extension of 3-dimensional Abelian Lie
 algebra $L^0$. For any $\xi,~\eta\in L^0$,
 $$
 [\xi,\eta]=\pi(\xi,\eta)=X\wedge \frac{\partial}{\partial x^4}(\xi,\eta)=0,
 $$
 since $\langle \frac{\partial}{\partial x^4},\xi\rangle=\langle \frac{\partial}{\partial x^4},\eta\rangle=0$.
  Furthermore we have $\langle e^4,[\xi,e_4]\rangle=0$ for any $\xi\in \frak g$ since $D(\pi)=0$
and this implies $[\xi,e_4]\in L^0$ and follows
 that $L^0$ is the Abelian idea of Lie algebra $\frak g$.

 As in the Case (2), the three linear independent lines $L_1,~L_2,~L_3$
 expand a 3-dimensional subspace $H$ and $H^0=ker(H)\in \frak g $ is a
 1-dimensional subspace of $\frak g$. Let $H^0=\mathbb{R}e_4$ and choose a
 3-dimensional subspace $E$ of $\frak g $ such that $\frak
 g=H^0\oplus E $. We will show $e_4$ is the center of Lie algebra $\frak
 g$ and follows that $\frak g$ is the central extension of Lie
 algebra $E$(mod $e_4$). Let $L_1=\mathbb{R}e^1,~ L_2=\mathbb{R}e^2,~ L_3=\mathbb{R}e^3$, so
$$
 \pi=\sum^4_{k=1}(C_{12}^kx^k \frac{\partial}{\partial x^1}\wedge \frac{\partial}{\partial x^2}+c_{13}^kx^k\frac{\partial}{\partial x^1}
\wedge \frac{\partial}{\partial
x^3}+C_{23}^kx^k\frac{\partial}{\partial x^2}\wedge
\frac{\partial}{\partial x^3}),
$$
 for some constants $C_{12}^k,~C_{13}^k,~C_{23}^k$. Obviously for any
 $\xi\in E$,
 $$[\xi,e_4]=\pi(\xi,e_4)=0,$$ so $e_4$ is the
 center of Lie algebra $\frak g$.

 As in the Case (3), choose a line $L\in P_l$ and with the same
 method in the Case (1), we can get the conclusion that it is the
 D-extension of some  3-dimensional Abelian Lie algebra. In fact, it is
 a particular case of case (1).

 When $D(\pi)=k\neq 0$, $k$ is the modular character of Lie algebra
 $\frak g$ defined as a vector in $\frak g^\ast$ such that
 $$
 \langle k,\xi\rangle=tr(ad (\xi)), ~\forall~\xi\in \frak g.
 $$
 Let $\frak h=\ker k$, $\forall~\xi\in \frak h,~tr(ad(\xi))=0.$  ~$\frak h^*=(ker k)^\ast=\frak g^\ast/\mathbb{R}k \subset \frak g^\ast$,
 the Lie-Poisson structure on $\frak h^*$ is just the reduction of $\Lambda_{\frak
 g}$ on $\frak h^*$. ~$ \forall
 ~\eta\in \frak g$, $$tr(ad([\xi,\eta]))=tr([ad(\xi),ad(\eta)])=0.$$ So $\frak h$ is the
 idea of the Lie algebra $\frak g$, and Lie algebra $\frak g$ is {\em
 D-extension} of Lie algebra $\frak h$.

 It is evident that the Abelian Lie algebra and Lie algebra $E$(mod $e_4$)
 is unimodular. The Lie-Poisson structure associated with $\frak h$ is just
 the reduction of $\Lambda_{\frak g}$ on $\frak h^*$, so $\frak h$ is unimodular. This
 complete the proof of the proposition.\qed

 \begin{rmk}
 ~If the modular character $k$ of the {\rm 4}-dimension Lie algebra is $(0,0,0,1)^T\in\mathbb R^4$, we should consider the D-extension of which the trace of the
 derivation $D$ is not zero, since the modular character of $\frak
 g_D$ is $(0,0,0,tr(D))^T$.
 \end{rmk}

So before we study the classification of linear Poisson structures
on $\mathbb R^4$, we should firstly study the classification of
linear Poisson structures on $\mathbb R^3$ which is done in
\cite{LLS} and the cohomology groups of 3-dimensional Lie algebras
with coefficients in trivial representation and adjoint
representation. This is the following two sections.

\section{The Classification of Lie-Poisson Structures on $\mathbb R^3$}

In this section, we list some useful result  in \cite{LLS} about the
classification of Lie-Poisson structures on $\mathbb R^3$.

In \cite{LLS}, it is pointed out that  Lie-Poisson structures
$\pi_{\frak g}$ on $\mathbb{R}^3$ are in one-to-one correspondence
with {\bf{compatible
        pair}} $(k,f)$ and denoted by $\pi_{k,f}$, where $k$ is the modular vector and $f$ is a
quadratic function, such that $\hat{k}f = 0$ and
\begin{equation}\label{eq pi k f}
 \pi_{\frak g}=\pi_{k,f}=\frac{1}{2}\widehat{I}\wedge\widehat{k}+\pi_{f}=
\frac{1}{2}\widehat{I}\wedge\widehat{k}+\frac{\partial f}{\partial
x}\py \wedge \pz+ \frac{\partial f}{\partial y}\pz \wedge \px+
\frac{\partial f}{\partial z}\px \wedge \py.
\end{equation}

\begin{thm}\label{Th. Isomor}{\rm\cite{LLS}}
            Let $\pi_{1}$ and $\pi_{2}$ be two linear Poisson structures on $\mathbb{R}^3$ determined
            by the compatible pairs
            $( k_{1},f_{1} )$ and $( k_{2},f_{2} )$ respectively. Then $\pi_{1}$ is
            isomorphic to $\pi_{2}$ if and only if there is a $T \in GL(3)$ such that
$$k_{2}=Tk_{1},\quad
                            f_{2}=det(T)f_{1}\circ T^{-1}.$$
\end{thm}

 \begin{cor}\label{cor autoDer} With notations above, consider the automorphism group of $\pi_\frkg$
 and derivation of $\frkg$,
 one has
            \begin{equation}
                Aut(\pi_\frkg) = \{T |~~ T \in GL(3),~~ Tk = k,~~f\circ T =
                det(T)f\},
            \end{equation}
 and
 \begin{equation}
                Der(\frkg) \cong \{D|~~ D \in {\frak gl(3)},~~ Dk = 0 ,~~  \hat{D}f =
                (trD)f\}.
 \end{equation}

        \end{cor}

\begin{thm}\label{Th classify}{\rm\cite{LLS}}
            Any Lie-Poisson structure $\pi_{k, f}$ on $\mathbb{R}^3$ is isomorphic to one of
            the following standard forms:

  \setlength{\parindent}{30mm} {\rm(A)}.~ $k=0$ (unimodular case),\qquad\qquad\quad {\rm(B)}.~$k=(0,0,1)^{T}$, i.e., $\hat{k}= \frac{\partial}{\partial z}$,
                  \begin{align*}
         (1).~f&=0         &     (7). ~ f&=0\\
         (2).~f&=x^2+y^2+z^2 &   (8).~f&=a(x^2+y^2)\\
          (3).~ f&=x^2+y^2-z^2 & (9).~f&=a(x^2-y^2)\\
           (4).~ f&=x^2+y^2 &    (10).~f&=x^2\\
         (5).~  f&=x^2-y^2 &\\
             (6).~ f&=x^2&
\end{align*}
where $a>0 $ is a constant.
 \end{thm}

\begin{thm}\label{Th Auto}{\rm\cite{LLS}}
        Let $G_i, ~ i = 1, \cdots, 10$ denote the automorphism groups of the Lie-Poisson structures which
        corresponds to Case (i) in Theorem \ref{Th classify}.
         Then we have
               \begin{itemize}
            \item[] $G_1 = GL(3)$,
            \item[] $G_2 = SO(3)$,
            \item[] $G_3 = SO(2,1)$,
            \item[] $G_4 = \{\left(
                                        \begin{array}
                                        {cc}
                                        \lambda T   &   0      \\
                                        \xi      &   det T

                                        \end{array}
                                \right)| ~~T \in
                                O(2), ~~\lambda\neq0, ~~\xi \in \mathbb{R}^2\}$,
            \item[] $G_5 = \{\left(
                                        \begin{array}
                                        {ccc}
                                        \alpha      &   \beta        &   0 \\
                                        \beta       &   \alpha       &   0 \\
                                        \gamma      &   \delta       &  1

                                    \end{array}
                                \right)or
                                \left(
                                        \begin{array}
                                        {ccc}
                                        -\alpha      &   \beta        &   0 \\
                                        -\beta       &   \alpha       &   0 \\
                                        \gamma      &   \delta       &  -1

                                    \end{array}
                                \right) |~~\alpha,~\beta,~\gamma,~\delta \in \mathbb{R},
                                 ~~\alpha^2\neq\beta^2\}$,
            \item[] $G_6 = \{\left(
                                        \begin{array}
                                        {cc}
                                        a           &   0     \\
                                        \xi       &   A

                                    \end{array}
                                \right)|~~A \in GL(2),
                                 ~~det(A) = a\neq0, ~~\xi \in \mathbb{R}^2\}$,
            \item[]  $G_7 = \{\left(
                                        \begin{array}
                                        {cc}
                                        A       &      0 \\
                                        \xi &      1

                                    \end{array}
                                \right)|~~A\in GL(2), ~~\xi \in \mathbb{R}^2\}$,
            \item[]  $G_8 = \{\left(
                                        \begin{array}
                                        {cc}
                                        \lambda T   &   0      \\
                                        \xi       &   1

                                        \end{array}
                                \right)| ~~T \in SO(2), ~~\lambda \neq 0, ~~\xi \in \mathbb{R}^2\}$,
            \item[]  $G_9 = \{\left(
                                        \begin{array}
                                        {ccc}
                                        \alpha      &   \beta        &   0 \\
                                        \beta       &   \alpha       &   0 \\
                                        \gamma      &   \delta       &   1

                                    \end{array}
                                \right) |~~\alpha,~\beta,~\gamma,~\delta \in \mathbb{R},~~ \alpha^2\neq\beta^2\}$,
            \item[]  $G_{10} = \{\left(
                                        \begin{array}
                                        {ccc}
                                        \alpha      &   0            &   0 \\
                                        \beta       &   \alpha       &   0 \\
                                        \gamma      &   \delta       &   1

                                    \end{array}
                                \right) |~~\alpha,~\beta,~\gamma,~\delta \in \mathbb{R}, ~~\alpha\neq
                                0\}$.
        \end{itemize}
    \end{thm}

\section{The Cohomology Groups of 3-dimensional Lie Algebras}

%As another application of the decomposition of  Lie-Poisson structures, we give a picture of  cohomology groups
%of 3-dimensional Lie algebras which will be used later.

In  this section, $\frak g $ will be a 3-dimensional Lie algebra
with modular character $k$ and $\frak g^*$ is its dual space with
Lie-Poisson structure $\pi_{k,f}$, described as in (\ref{eq pi k
f}), where $f$ is decided by  symmetric matrix $A$.
%Denote by $\omega$ a Chevalley-Eilenberg 2-chain of $\frak g$ and by $D$ a 1-chain corresponding to adjoint representation.

\begin{thm}\label{thm dw}
For any  $\eta \in \frak g^*$, $\omega\in \frak g^*\wedge \frak
g^*$,
%where $\frak g$ is a
 %3-dimensional Lie algebra and $\frak g^*$ is its dual space, $k$ is the modular
%character of $\frak g$, $A$ is the symmetric matrix decided by $f$,
%where the pair$(k,f)$ is in correspondence with the Lie-Poisson structure$\pi_{\frak g}$ on $\frakg^*$,
$\delta$ is the coboundary operator of Lie algebra cohomology
with coefficients in its trivial representation,  then
\begin{equation}\label{d eta}
\delta\eta=\frac{1}{2}k\wedge\eta-2\Phi \circ A\eta,
\end{equation}
\begin{equation}\label{dw}
\delta\omega=k\wedge\omega.
\end{equation}

\end{thm}
\pf By Theorem \ref{Th. decom}, $\pi =
\frac{1}{2}\hat{I}\wedge\hat{k} + \Lambda_{{\frak g}}$, and
$$
~\Lambda_{{\frak g}}=\pi_f=\frac{\partial f}{\partial x_1}
\frac{\partial}{\partial x_2}\wedge\frac{\partial}{\partial
x_3}+\frac{\partial f}{\partial x_2}\frac{\partial}{\partial
x_3}\wedge\frac{\partial}{\partial x_1}+\frac{\partial f}{\partial
x_3}\frac{\partial}{\partial x_1}\wedge\frac{\partial}{\partial
x_2}.
$$
So we have
$$
[e_1,e_2]=\frac{1}{2}(\langle k,e_2\rangle e_1-\langle k,e_1\rangle
e_2)+2Ae^3,$$$$ [e_2,e_3]=\frac{1}{2}(\langle k,e_3\rangle
e_2-\langle k,e_2\rangle e_3)+2Ae^1,$$$$
[e_3,e_1]=\frac{1}{2}(\langle k,e_1\rangle e_3-\langle k,e_3\rangle
e_1)+2Ae^2,
$$
and \begin{eqnarray}\nonumber \langle \delta\eta,e_1\wedge e_2\rangle &=&-\langle\eta,[e_1,e_2]\rangle\\
\nonumber &=&-\langle\eta,\frac{1}{2}(\langle k,e_2\rangle
e_1-\langle k,e_1\rangle
e_2)+2Ae^3\rangle\\
\nonumber &=&\frac{1}{2}\langle
k,e_1\rangle\langle\eta,e_2\rangle-\frac{1}{2}\langle
k,e_2\rangle\langle\eta,e_1\rangle-2\langle A\eta,\Phi(e_1\wedge e_2)\rangle\\
\nonumber &=&\langle\frac{1}{2}k\wedge \eta, e_1\wedge
e_2\rangle-2\langle\Phi\circ A\eta,e_1\wedge e_2\rangle.
\end{eqnarray}
After the same computation of $\langle\delta\eta,e_2\wedge
e_3\rangle$ and $\langle\delta\eta,e_3\wedge e_1\rangle$, we have
the conclusion that
$$\delta\eta=\frac{1}{2}k\wedge\eta-2\Phi \circ A\eta.$$
 For
$\eta_1,~\eta_2\in \frak g^*$, by Leibnitz rule, we have
\begin{eqnarray}
\nonumber\delta(\eta_1\wedge\eta_2)&=&\delta\eta_1\wedge\eta_2-\eta_1\wedge\delta\eta_2\\
\nonumber&=&(\frac{1}{2}k\wedge\eta_1-2\Phi \circ
A\eta_1)\wedge\eta_2-\eta_1\wedge(\frac{1}{2}k\wedge\eta_2-2\Phi
\circ A\eta_2)\\
\nonumber&=&k\wedge\eta_1\wedge\eta_2+2(\eta_1\wedge(\Phi\circ
A\eta_2)-(\Phi \circ
A\eta_1)\wedge\eta_2)\\\nonumber&=&k\wedge\eta_1\wedge\eta_2.
\end{eqnarray}
Thus for any $\omega \in \frak g^*\wedge \frak g^*$,
$\delta\omega=k\wedge\omega$. \qed

For Lie algebras $\frak g$ listed in Theorem \ref{Th classify},
 we  give a detail description of the corresponding
2-cocycles, denote by $\mathcal C^2(\frak g)\subset \frak g^*\wedge
\frak g^*$ and exact 2-cocycles, denote by $\mathcal B^2(\frak
g)\subset \frak g^*\wedge \frak g^*$, which will be used when
consider the classification of linear Poisson structures on $\mathbb
R^4$. Furthermore we give cohomology groups $H^1(\frak g),~H^2(\frak
g),~H^3(\frak g)$ with coefficients in trivial representation by the
way.

\begin{cor}\label{w exact}\hspace{\fill}

(A).~$k=0$(unimodular case)

\begin{tabular}{|l|l|l|l|l|l|}\hline
\multicolumn{1}{|c|}{quadratic function
$f$}&\multicolumn{1}{|c|}{$\mathcal{C}^2(\frak g)$}
&\multicolumn{1}{|c|}{$\mathcal{B}^2(\frak g)$}&
\multicolumn{1}{|c|}{$H^1({\frak
g})$}& \multicolumn{1}{|c|}{$H^2({\frak g})$}& \multicolumn{1}{|c|}{$H^3({\frak g})$}\\
\hline

\ZZ{-4pt}{18pt}$(1).~0$ & $\forall ~\omega$  & 0 & $\mathbb R^3$&
$\mathbb R^3$& $\mathbb R$\\\hline

\ZZ{-4pt}{18pt}$(2).~x^2+y^2+z^2$  & $\forall ~\omega$  & $\forall
~\omega$ & 0& 0& $\mathbb R$\\\hline

\ZZ{-4pt}{18pt}$(3).~x^2+y^2-z^2$ & $\forall ~\omega$  & $\forall
~\omega$ & 0& 0& $\mathbb R$
\\\hline

\ZZ{-4pt}{18pt}$(4).~x^2+y^2$  & $\forall ~\omega$ & $\alpha
dy\wedge dz+\beta dz\wedge dx$ & $\mathbb R$& $\mathbb R$& $\mathbb
R$\\\hline

\ZZ{-4pt}{18pt}$(5).~x^2-y^2$ &$\forall ~\omega$ & $\alpha dy\wedge
dz+\beta dz\wedge dx$ & $\mathbb R$& $\mathbb R$& $\mathbb
R$\\\hline

\ZZ{-4pt}{18pt}$(6).~x^2$ & $\forall ~\omega$ & $\alpha dy\wedge
dz$& $\mathbb R^2$& $\mathbb R^2$& $\mathbb R$\\\hline
\end{tabular}

(B).~$k=(0,0,1)^{T}$, i.e., $\hat{k}= \frac{\partial}{\partial z}$

\begin{tabular}{|l|l|l|l|l|l|}\hline
\multicolumn{1}{|c|}{quadratic function
$f$}&\multicolumn{1}{|c|}{$\mathcal{C}^2(\frak
g)$}&\multicolumn{1}{|c|}{$\mathcal{B}^2(\frak
g)$}&\multicolumn{1}{|c|}{ $H^1({\frak g})$}&\multicolumn{1}{|c|}{
$H^2({\frak g})$}& \multicolumn{1}{|c|}{$H^3({\frak g})$}\\\hline

\ZZ{-4pt}{18pt}$(7).~0 $& $\alpha dy\wedge dz+\beta dz\wedge dx$
&$\alpha dy\wedge dz+\beta dz\wedge dx$ & $\mathbb R$ & 0& 0
\\\hline

\ZZ{-4pt}{18pt}$(8).~a(x^2+y^2)$  & $\alpha dy\wedge dz+\beta
dz\wedge dx$ & $\alpha dy\wedge dz+\beta dz\wedge dx$ & $\mathbb R$
& 0& 0
\\\hline

\ZZ{-4pt}{18pt}$(9)_1.~a(x^2-y^2),~a\neq\frac{1}{4}$ &$\alpha
dy\wedge dz+\beta dz\wedge dx$ & $\alpha dy\wedge dz+\beta dz\wedge
dx$ & $\mathbb R$ & 0& 0\\\hline

\ZZ{-4pt}{18pt}$(9)_2.~\frac{1}{4}(x^2-y^2)$&$\alpha dy\wedge
dz+\beta dz\wedge dx$ &$\alpha (dy\wedge dz- dz\wedge dx)$&$\mathbb
R^2$&$\mathbb R$&0
\\\hline

\ZZ{-4pt}{18pt}$(10).~x^2$ &$\alpha dy\wedge dz+\beta dz\wedge dx$ &
$\alpha dy\wedge dz+\beta dz\wedge dx$& $\mathbb R$ & 0& 0\\\hline
\end{tabular}

where $\alpha,~\beta$ are arbitrary constants.
\end{cor}

\pf Firstly consider $H^3$, by Theorem \ref{thm dw}, for $\forall~
\omega\in\frak g^*\wedge\frak g^*$, $\delta\omega=k\wedge\omega$. So
in the unimodular case $k=0$, this implies $\delta\omega=0$. So
$\theta\in \wedge^3 \frak g^*$ is exact if and only if $\theta=0$,
and this implies $H^3=\mathbb{R}$; If $k=(0,0,1)$, $\forall
~\theta\in \wedge^3 \frak g^*$ is exact, so $H^3=0$.

Next consider $H^1$, there is no exact chain and $\eta\in \frak g^*$
is closed if and if
$$
\frac{1}{2}k\wedge\eta-2\Phi \circ A\eta=0
$$ by
Theorem \ref{thm dw}. If $k=0$, this equals $A\eta=0$, so the
dimension of the first cohomology group is $3-order(A)$ and we have
the conclusion listed above. If $k=(0,0,1)$, it is a little
complicated but  straightforward however and we leave it to the
interesting reader.

At last we consider $H^2$, in the case that $k=0$, all of 2-chains
are 2-cocycles by Theorem \ref{thm dw}; In the case that
$k=(0,0,1)$, $\omega\in \wedge^2 \frak g^*$ is closed if and only if
$\omega$ has the form $\omega=\alpha dy\wedge dz+\beta dz\wedge dx$
by Theorem \ref{thm dw}, where $\alpha, \beta\in \mathbb{R}$.
$\omega\in \wedge^2 \frak g^*$ is exact if and only if
$$
\omega =\frac{1}{2}k\wedge\eta-2\Phi \circ A\eta
$$
for some $\eta
\in \frak g^* $ by Theorem \ref{thm dw}, the conclusion is
straightforward.\qed

\begin{rmk}
In fact, the first cohomology group $H^1(\frak g)$ relate to the
dimension of derived algebra $[\frak g,\frak g]$ of 3-dimensional
Lie algebra $\frak g$, more precise, $dim(H^1(\frak g))=3-dim([\frak
g,\frak g])$. In  the unimodular case, $dim([\frak g,\frak
g])=order(A)$, however it is not true if the modular vector is not
zero, see the case that $f=\frac{1}{4}(x^2-y^2)$.
\end{rmk}

The next theorem and corollary give some description of the
derivation which will be used in the next section when we consider
the extension of a Lie algebra by a derivation.

\begin{thm}\label{thm D }
Let $\frak g $ be one of the Lie algebras listed in Theorem \ref{Th
classify}. $D$ is a derivation of Lie algebra $\frak g$ if and only
if
\begin{equation}
D^*k=0,\quad DA+AD^*=tr(D)A.
\end{equation}
And $D$ is an inner derivation if and only if there exists a
skew-symmetric transformation $B:\frak g\rightarrow \frak g^*$ such
that
\begin{equation}
D=(2A+\frac{1}{2}\overline{k})B,
\end{equation}
where $\overline{k}:\frak g^*\longrightarrow \frak g$ denotes the
skew-symmetric transformation
 decided by $\Phi^{-1}k$.
\end{thm}
\pf ~Since $D$ is a derivation, by Corollary \ref{cor autoDer} we
have $$D^*k=0,~\hat{D^*}f =tr(D)f,$$ and this implies
$$D^*k=0,~DA+AD^*=tr(D)A.$$
If $D$ is an inner derivation, there exists a $\xi\in \frak g$, such
that $D=ad(\xi)$ and $DX=ad(\xi)(X)=[\xi,X]$. Furthermore,
 \begin{eqnarray}
 \nonumber[\xi,X]&=&c.p.\{(\langle \xi,e^1\rangle\langle X,e^2\rangle-\langle\xi,e^2\rangle\langle X,e^1\rangle)[e_1,e_2]\}\\
  \nonumber &=&
  c.p.\{(\langle\xi,e^1\rangle\langle X,e^2\rangle-\langle\xi,e^2\rangle\langle X,e^1\rangle)(\frac{1}{2}(\langle k,e_2\rangle e_1-\langle k,e_1\rangle e_2)+2Ae^3)\}\\
  \nonumber &=&c.p.\{\frac{1}{2}\Phi^{-1}(k\wedge
  (\langle\xi,e^1\rangle\langle X,e^2\rangle-\langle\xi,e^2\rangle\langle X,e^1\rangle)e^3)\}+2A\Phi(\xi\wedge X) \\
  \nonumber &=&\frac{1}{2}\Phi^{-1}(k\wedge(\Phi(\xi \wedge X)))+ 2A\Phi(\xi\wedge
  X)
  \end{eqnarray}
where $c.p.\{\}$ means the cyclic permutations of $e^1,~ e^2,~ e^3$.
So we have
\begin{eqnarray*}
D(\cdot)&=&\frac{1}{2}\Phi^{-1}(k\wedge(\Phi(\xi)(\cdot)))+2A\Phi(\xi)(\cdot)\\
&=&(2A+\frac{1}{2}\overline{k})B(\cdot),
\end{eqnarray*}
where $\overline{B}:\frak g\rightarrow \frak g^*$ is the
skew-symmetric transformation   decided by $\Phi\xi$ and
$\overline{k}:\frak g^*\rightarrow \frak g$ is the skew-symmetric
transformation  decided by $\Phi^{-1}k$. \qed

Now we can  give the form of the derivation and the inner derivation
via the standard form of Lie algebras listed  in Theorem \ref{Th
classify}. However inner derivation will be of  nothing use when we
consider the extension of a Lie algebra by a derivation.
\newpage

\begin{cor}\label{D exact}\hspace{\fill}

(A).~$k=0$(unimodular case)

\begin{tabular}{|l|l|l|l|}\hline
\multicolumn{1}{|c|}{quadratic function
$f$}&\multicolumn{1}{|c|}{derivation}&\multicolumn{1}{|c|}{inner
derivation}&\multicolumn{1}{|c|}{$H^1(ad;~\frak g)$}\\\hline

\ZZ{-4pt}{18pt} $ (1).~0$ &  $\forall~D\in \frak gl(3) $ & 0 &
$\mathbb R^9$\\\hline

\ZZ{-4pt}{18pt} $(2).~x^2+y^2+z^2$  &  $D\in \frak o(3)$ & $D\in
\frak o(3)$& 0\\\hline

\ZZ{-4pt}{18pt} $(3).~x^2+y^2-z^2$  &  $D\in \frak o(2,1)$ & $D\in
\frak o(2,1)$&0\\\hline

\ZZ{-20pt}{45pt} $(4).~x^2+y^2$  & $\left(\begin{array}{ccc}
\alpha& \beta& \gamma\\
-\beta& \alpha& \delta\\
0& 0& 0
\end{array}
\right)$ & $\alpha=0$& $\mathbb R$\\\hline

\ZZ{-20pt}{45pt} $(5).~x^2-y^2$ & $\left(\begin{array}{ccc}
\alpha& \beta& \gamma\\
\beta& \alpha& \delta\\
0& 0& 0
\end{array}
\right)$& $\alpha=0$& $\mathbb R$\\\hline

\ZZ{-20pt}{45pt} $(6).~x^2$ & $\left(\begin{array}{cc}
tr(\overline{D})& \xi \\
0& \overline{D}
\end{array}\right)$, $\xi \in \mathbb R^2$, $\overline{D}\in
\frak gl(2)$ & $\overline{D}=0$& $\mathbb R^4$\\\hline
\end{tabular}

(B).~$k=(0,0,1)^{T}$, i.e., $\hat{k}= \frac{\partial}{\partial z}$

\begin{tabular}{|l|l|l|l|}\hline
\multicolumn{1}{|c|}{quadratic function
$f$}&\multicolumn{1}{|c|}{derivation}&\multicolumn{1}{|c|}{inner
derivation} & \multicolumn{1}{|c|}{$H^1(ad;~\frak g)$}\\\hline

\ZZ{-20pt}{45pt} $(7).~0$ &
$\left(\begin{array}{ccc}\alpha&\beta&\gamma\\\delta &\epsilon
&\varepsilon\\0&0&0\end{array}\right)$ &
$\left(\begin{array}{ccc}\epsilon&o&\gamma\\0&\epsilon
&\varepsilon\\0&0&0\end{array}\right)$&$\mathbb R^3$\\\hline

\ZZ{-20pt}{45pt} $(8).~a(x^2+y^2)$  &
$\left(\begin{array}{ccc}\alpha&\beta&\gamma\\-
                   \beta&\alpha&\delta\\0&0&0\end{array}\right)$&
$\left(\begin{array}{ccc}\alpha&4a\alpha&\gamma\\-4a\alpha&\alpha&\delta\\0&0&0\end{array}\right)$&
$\mathbb R$ \\\hline

\ZZ{-20pt}{45pt} $(9)_1.~a(x^2-y^2),~a\neq\frac{1}{4}$ &
$\left(\begin{array}{ccc}\alpha&\beta&\gamma\\
                   \beta&\alpha&\delta\\0&0&0\end{array}\right)$&
$\left(\begin{array}{ccc}\alpha&4a\alpha&\gamma\\4a\alpha&\alpha&\delta\\0&0&0\end{array}\right)$&
$\mathbb R$\\\hline

\ZZ{-20pt}{45pt} $(9)_2.~\frac{1}{4}(x^2-y^2),$ &
$\left(\begin{array}{ccc}\alpha&\beta&\gamma\\
                   \beta&\alpha&\delta\\0&0&0\end{array}\right)$&
$\left(\begin{array}{ccc}\alpha&\alpha&\gamma\\\alpha&\alpha&\gamma\\0&0&0\end{array}\right)$&
$\mathbb R^2$\\\hline

\ZZ{-20pt}{45pt} $(10).~x^2$ &
$\left(\begin{array}{ccc}\alpha&\beta&\gamma\\0&\alpha&\delta\\0&0&0\end{array}\right)$
& $
\left(\begin{array}{ccc}\alpha&4\alpha&\gamma\\0&\alpha&\delta\\0&0&0\end{array}\right)$&$\mathbb
R$\\\hline
\end{tabular}

where $a>0,~\alpha,\beta,\gamma,\delta,\epsilon,\varepsilon\in
\mathbb R$ are arbitrary constants.

\end{cor}
\pf The proof is almost straightforward and we only give the proof
of the Case (2). Any derivation $D$ satisfies
         $DA+AD^*=tr(D^*)A$.  Multiply $A^{-1}$ at right, we have $D=-AD^*A^{-1}+tr(D^*)I$, and this implies
                              $$tr(D)=-tr(D)+3tr(D)\Longrightarrow tr(D)=0.$$ Then $D=-AD^*A^{-1}$, this implies $D$
                              is skew-symmetric. To prove
                              $D$ is an inner derivation, we only need to say $-D^*A^{-1}$ is skew-symmetric.
                              Notice that $DA+AD^*=0$ means $D^*A^{-1}+A^{-1}D=0$, this is just
                              the condition that $-D^*A^{-1}$ is
                              skew-symmetric. \qed

\section{The Classification of  Linear Poisson Structures  on $\mathbb R^4$}

With above preparations,  we can see that the procedure of
classification of Lie-Poisson structures on $\mathbb{R}^4$ can be
split into four steps by Theorem \ref{Th classify}, Corollary \ref{w
exact} and \ref{D exact} as follows.
\begin{itemize}
                    \item[(1)] Take a standard form of   Lie-Poisson
                    structure $\pi_{{\frak g}}$
                    from the list in Theorem \ref{Th classify}(A).
                     \item[(2)]  Consider all the isomorphism
                     classes
                     of the extension, central extension and
                     D-extension  of the corresponding
                     Lie algebra.
                      \item[(3)] Consider the isomorphism of the
                      extension of two different Lie algebras.
                     \item[(4)] The Lie-Poisson structures
                     corresponding to the above Lie algebras give the
                     classification of linear Poisson structures on
                     $\mathbb R^4$
\end{itemize}

In fact Step 2 above is the most difficult and complicated one. If
this is done,  Step 3 is almost straightforward however. The
following proposition  give a detail description of central
extensions and D-extensions of  Lie algebras listed in Theorem
\ref{Th classify} and of which the first part will be used in the
classification of 4-dimensional Lie-Poisson structures and the whole
proposition will be of great importance when we consider {\em affine
Poisson structures}.

\begin{pro}\label{thm iso}%%%%%%%%%%%%%%%%%%%%%%%%%%%%%%%%定理3.6
 With the same notations above, Let $\frak g$ be one of the
 3-dimensional Lie algebras listed in Theorem \ref{Th
classify},   any of its extension, central extension which is
decided by 2-cocycles $\omega$ in trivial representation and
D-extension which is decided by derivations, is isomorphic to one of
the following forms:

(A).~$k=0$(unimodular case)\\[3pt]
\begin{tabular}{|l|l|l|}\hline
\multicolumn{1}{|c|}{quadratic function
$f$}&\multicolumn{1}{|c|}{2-cocycles ~$
\omega$}&\multicolumn{1}{|c|}{Derivation D}\\\hline

\ZZ{-40pt}{85pt} $(1).~0$ & 0 and $ dx\wedge dy$ &
\begin{minipage}{95mm}0,
$\left(\begin{array}{ccc}0&1&0\\0&0&1\\0&0&0\end{array}\right),~\left(\begin{array}{ccc}0&1&0\\0&0&0\\0&0&0\end{array}\right),~
~\left(\begin{array}{ccc}\alpha&1&0\\-1&\alpha&0\\0&0&\beta\end{array}\right)$\\
$\left(\begin{array}{ccc}1&0&0\\0&0&1\\0&0&0\end{array}\right),\left(\begin{array}{ccc}1&1&0\\0&1&0\\0&0&\beta\end{array}\right),
\left(\begin{array}{ccc}1&0&0\\0&\alpha&0\\0&0&\beta\end{array}\right),\left(\begin{array}{ccc}1&1&0\\0&1&1\\0&0&1\end{array}\right)$
\end{minipage}\\\hline

\ZZ{-4pt}{18pt} $(2).~x^2+y^2+z^2$  &  0 & 0\\\hline

\ZZ{-4pt}{18pt} $(3).~x^2+y^2-z^2$  &  0 & 0\\\hline

\ZZ{-20pt}{45pt} $(4).~x^2+y^2$  & 0 and $ dx\wedge dy$ & 0 or
$\left(\begin{array}{ccc}1&0&0\\0&1&0\\0&0&0\end{array}\right)$\\\hline

\ZZ{-20pt}{45pt} $(5).~x^2-y^2$ &0 and $ dx\wedge dy$ & 0 or
$\left(\begin{array}{ccc}1&0&0\\0&1&0\\0&0&0\end{array}\right)$\\\hline

\ZZ{-40pt}{85pt} $(6).~x^2$ & 0 and $ dx\wedge dy$
&\begin{minipage}{85mm}
0,~$\left(\begin{array}{ccc}0&0&0\\0&0&1\\0&0&0\end{array}\right),~\left(\begin{array}{ccc}0&0&0\\0&1&0\\0&0&-1\end{array}\right),~\left(\begin{array}{ccc}0&0&0\\0&0&1\\0&-1&0\end{array}\right)
$,\\
$\left(\begin{array}{ccc}2&0&0\\0&1&\alpha\\0&-\alpha&1\end{array}\right),~\left(\begin{array}{ccc}2&0&0\\0&1&1\\0&0&1\end{array}\right),~\left(\begin{array}{ccc}1&0&0\\0&\alpha&0\\0&0&1-\alpha\end{array}\right)$
\end{minipage}\\\hline
\end{tabular}

\newpage(B).~$k=(0,0,1)^{T}$, i.e., $\hat{k}= \frac{\partial}{\partial
z}$\\[3pt]
\begin{tabular}{|l|l|l|}\hline
\multicolumn{1}{|c|}{quadratic function
$f$}&\multicolumn{1}{|c|}{2-cocycles
$\omega$}&\multicolumn{1}{|c|}{Derivation D}\\\hline

\ZZ{-40pt}{85pt} $(7).~0$ & 0
 &\begin{minipage}{85mm} 0,
 $\left(\begin{array}{ccc}0&1&0\\0&0&0\\0&0&0\end{array}\right),~\left(\begin{array}{ccc}1&1&0\\0&1&0\\0&0&0\end{array}\right),~
\left(\begin{array}{ccc}1&0&0\\0&\alpha&0\\0&0&0\end{array}\right),$\\$\left(\begin{array}{ccc}0&\alpha&0\\-\alpha&0&0\\0&0&0\end{array}\right),~
\left(\begin{array}{ccc}1&\alpha&0\\-\alpha&1&0\\0&0&0\end{array}\right)$
\end{minipage}\\\hline

\ZZ{-20pt}{45pt} $(8).~a(x^2+y^2)$  & 0& 0 and
$\left(\begin{array}{cc}I_{2\times
2}&0\\0&0\end{array}\right)$\\\hline

\ZZ{-20pt}{45pt} $(9)_1.~a(x^2-y^2),~a\neq\frac{1}{4}$ & 0 &0 and
$\left(\begin{array}{cc}I_{2\times
2}&0\\0&0\end{array}\right)$\\\hline

\ZZ{-20pt}{45pt} $(9)_2.~\frac{1}{4}(x^2-y^2),$ & 0 and $dx\wedge
dz$ & 0,~$\left(\begin{array}{cc}I_{2\times
2}&0\\0&0\end{array}\right),~\left(\begin{array}{ccc}0&1&0\\0&0&0\\0&0&0\end{array}\right)$\\\hline

\ZZ{-20pt}{45pt} $(10).~x^2$ & 0 &0 and
$\left(\begin{array}{ccc}0&0&1\\0&0&0\\0&0&0\end{array}\right)$\\\hline
\end{tabular}\\[5pt]
where $a> 0,~\alpha,~\beta\in\mathbb R $ are constants.
\end{pro}

\pf Throughout the whole proof, 2-chain $\omega=\alpha dy\wedge
dz+\beta dz\wedge dx+\gamma dx\wedge dy$ will be denoted by
$\omega=(\alpha,\beta,\gamma)$.

(1).~~ Assume $\omega_1=(\alpha_1,\beta_1,\gamma_1)$,~
$\omega_2=(\alpha_2,\beta_2,\gamma_2)$, where
$\alpha_1,~\beta_1,~\gamma_1$ are not zero at the same time and so
are $\alpha_2,~\beta_2,~\gamma_2$, and $A$  is a matrix of which the
adjoint matrix
    $\widetilde{A}$ satisfies $\widetilde{A}\omega_2=\omega_1$, then
    $\overline{A}=\left(\begin{array}{cc}A&0\\0&1\end{array}\right)$
    is the isomorphism from $\frak g_{\omega_1}$ to $\frak
    g_{\omega_2}$. In particular, we can choose  $\omega=dx\wedge dy$ as the standard form.

    (2) and (3) follows from  Corollary \ref{w exact} and Corollary \ref{D exact}.

    (4) and (5). By Corollary  \ref{w exact} and the fact that if  $\omega_1-\omega_2$ is exact, $\frak g_{\omega_1}$ is
          isomorphic to  $\frak g_{\omega_2}$, we only need to
          consider the case $\omega=\gamma dx\wedge dy$, where $\gamma\neq 0$.
          Let $\omega_1=\gamma_1dx\wedge dy$ and $\omega_2=\gamma_2dx\wedge
          dy$,~ $\gamma_1\neq0,~ \gamma_2\neq0$,~  $\frak g_{\omega_1}$ is isomorphic to  $\frak
          g_{\omega_2}$ is obvious since
          $A=\left(\begin{array}{cc}I_{3\times3}&0\\0&\frac{\gamma_2}{\gamma_1}\end{array}\right)$
          is the isomorphism and we can choose $\omega=dx\wedge dy$ as the standard form; By  Corollary \ref{D exact} and
           the fact that if $D_1-D_2$ is exact, $\frak
          g_{D_1}$ is isomorphic to  $\frak  g_{D_2}$,  we only need
          to consider the case
          $D=\left(\begin{array}{cc}   d\cdot I_{2\times 2}&0\\0&0
          \end{array}\right)$, where $d\neq0$.  $\frak  g_{D_1}$ is isomorphic to  $\frak
          g_{D_2}$ is obvious since  $D=\left(\begin{array}{cc}I_{3\times3}&0\\0&\frac{d_1}{d_2}\end{array}\right)$
           is the isomorphism. And we can choose  $D=\left(\begin{array}{cc}
                                                 I_{2\times 2}&0\\0&0
                                                  \end{array}\right)$
                                                  as the standard
                                                  form.

      (6). As in {\rm (4)} and {\rm (5)}, we only need to consider the case $\omega=
      (0,\beta,\gamma) $,  assume  $\omega_1=(0,\beta_1,\gamma_1) $ and   $\omega_2=(0,\beta_2, \gamma_2)$,
       Let $B\in GL(2)$ that satisfies $B\omega_2=\omega_1$, then
      $A=\left(\begin{array}{ccc}
      \sqrt[3]{det(B)}&0&0\\
      0&\frac{1}{\sqrt[3]{det(B)}}B&0\\
      0&0&1\end{array}\right)$   is the isomorphism from $\frak g_{\omega_1}$
      to $\frak g_{\omega_2}$. As for {\em D-extension }, we only
      need to consider the case  $D=\left(\begin{array}{cc}tr(\overline{D})& 0\\ 0&
      \overline{D}\end{array}\right)$ by Corollary \ref{D exact}, and $\frak g_{D_1}$
      is isomorphic to $\frak g_{D_2}$ if and only if   $\overline{D_1}$
      is similar to the matrix that is nonzero  multiples of the matrix similar
      to  $\overline{D_2}$. In fact, if $\overline{D_1}=d\overline{B}^{-1}\overline{D_2}\overline{B}$ for some $\overline{B}\in GL(2)$,
      then $B=\left(\begin{array}{ccc} det(\overline{B})& 0&0\\ 0&\overline{B}&0\\0&0&d\end{array}\right)$
      is the isomorphism from $\frak g_{D_1}$ to $\frak g_{D_2}$.

      The proof  of the triviality of central extensions in the Cases
(7),~(8),~(10) and in the Case $(9)_1$  are same because of
Corollary \ref{w exact}. In the Case $(9)_2$, if the 2-cocycle
$\omega$ is exact, from Corollary \ref{w exact}, we have
$\alpha=-\beta$. So we only need to consider the case
$\omega=\beta\frac{\partial}{\partial
z}\wedge\frac{\partial}{\partial x}$ and $\frak g_{\omega_1}$ is
isomorphic to  $\frak g_{\omega_2}$ is obvious since
$B=\left(\begin{array}{ccc}I_{3\times
3}&0\\0&\frac{\beta_2}{\beta_1}\end{array}\right)$ is the
isomorphism, where $\omega_1=\beta_1\frac{\partial}{\partial
z}\wedge\frac{\partial}{\partial x} $,~~
$\omega_2=\beta_2\frac{\partial}{\partial
z}\wedge\frac{\partial}{\partial x} $.

The proof of the determination of derivation $D$ in Case $(B)$ when
$k=(0,0,1)^T$ listed in Theorem \ref{Th classify} is almost the same
as the proof of the unimodular case and we leave it to the
interesting reader.\qed

Now combine Proposition  \ref{extension}, Theorem \ref{Th classify}
and \ref{thm iso}, we obtain the following theorem that gives the
classification of Lie-Poisson structures on $\mathbb R^4$.

\begin{thm}\label{thm main}
Any 4-dimensional Lie algebra is isomorphic to one of the following
forms which are  considered as the extension by a derivation of some
unimodular 3-dimensional Lie algebra listed in Theorem \ref{Th
classify}(A).

\begin{tabular}{|l|l|}\hline
\multicolumn{1}{|c|}{quadratic function
$f$}&\multicolumn{1}{|c|}{Derivation D}\\\hline

\ZZ{-40pt}{85pt} $(1).~0$ &  \begin{minipage}{75mm}0,~
$\left(\begin{array}{ccc}0&1&0\\0&0&1\\0&0&0\end{array}\right),~\left(\begin{array}{ccc}1&0&0\\0&0&1\\0&0&0\end{array}\right),~\left(\begin{array}{ccc}1&1&0\\0&1&1\\0&0&1\end{array}\right),~
$\\
$\left(\begin{array}{ccc}1&1&0\\0&1&0\\0&0&\alpha\end{array}\right),~\left(\begin{array}{ccc}1&0&0\\0&\beta&0\\0&0&\gamma\end{array}\right),~
~\left(\begin{array}{ccc}\delta&1&0\\-1&\delta&0\\0&0&\epsilon\end{array}\right)$
\end{minipage}\\\hline

\ZZ{-4pt}{18pt} $(2).~x^2+y^2+z^2$   & 0\\\hline

\ZZ{-4pt}{18pt} $(3).~x^2+y^2-z^2$   & 0\\\hline

\ZZ{-20pt}{45pt} $(4).~x^2+y^2$   & 0 and
$\left(\begin{array}{ccc}1&0&0\\0&1&0\\0&0&0\end{array}\right)$\\\hline

\ZZ{-20pt}{45pt} $(5).~x^2-y^2$
 & 0 and
$\left(\begin{array}{ccc}1&0&0\\0&1&0\\0&0&0\end{array}\right)$\\\hline

\ZZ{-40pt}{85pt} $(6).~x^2$  &\begin{minipage}{85mm}
0,~$\left(\begin{array}{ccc}0&0&0\\0&1&0\\0&0&-1\end{array}\right),~\left(\begin{array}{ccc}0&0&0\\0&0&1\\0&-1&0\end{array}\right),
~\left(\begin{array}{ccc}2&0&0\\0&1&1\\0&0&1\end{array}\right)
$,\\
$\left(\begin{array}{ccc}2&0&0\\0&1&\varepsilon\\0&-\varepsilon&1\end{array}\right),~\left(\begin{array}{ccc}1&0&0\\0&\zeta&0\\0&0&1-\zeta\end{array}\right)$
\end{minipage}\\\hline
\end{tabular}

where
$\alpha,~\beta,~\gamma,~\delta,~\epsilon,~\varepsilon,~\zeta\in
\mathbb R$ are arbitrary constants.
\end{thm}

\pf~ By Proposition \ref{extension}, we know any 4-dimensional Lie
algebra is the extension of some 3-dimensional unimodular Lie
algebra and Proposition \ref{thm iso} gives the classification of
the extension of a fixed Lie algebra. Then only thing left now is to
decide the isomorphism of the extension of different Lie algebras,
however this is easy to be done. This complete the proof.\qed

%\begin{rmk}%%%%%%%%%%%%%%%%%%%%%% 注
%In fact, all the central extensions above can be considered as some D-extension. For example, the central extension of Lie algebra
%$\frak g$ decided by $f=x_1^2+x_2^2$ have the standard form$[x_2,x_3]=2x_1, [x_1,x_3]=-2x_2, [x_1,x_2]=x_4$, and this Lie
%algebra can be considered as the D-extension of the Lie algebra expanded by $x_1,x_2,x_4$ and the corresponding derivation
%$D=\left(\begin{array}{ccc}0& 2&  0\\-2& 0& 0\\0& 0& 0\end{array}\right)$.\end{rmk}

 By Proposition  \ref{thm iso} and the statement at the beginning of this section, we classify all of the {\em
affine Poisson structures} on $ \mathbb R^3$.

\begin{thm} On $\mathbb R^3$, any {\em affine Poisson structure} is
isomorphic to one of the following forms:

(A).~$k=0$(unimodular case)
\qquad\qquad\qquad\qquad\quad(B).~$k=(0,0,1)^{T}$, i.e., $\hat{k}=
\frac{\partial}{\partial z}$\\
\begin{tabular}{|l|l|}\hline
\ZZ{-4pt}{18pt}quadratic function $f$&affine Poisson structure
\\\hline

\ZZ{-4pt}{18pt}$(1).~0$ & $\frac{\partial}{\partial
x}\wedge\frac{\partial}{\partial y} $
\\\hline

\ZZ{-6pt}{18pt}$(2),(3).~x^2+y^2\pm z^2$  &  $\pi_f$\\\hline

\ZZ{-6pt}{18pt}$(4).~x^2+y^2$  & $\frac{\partial}{\partial
x}\wedge\frac{\partial}{\partial y}+\pi_f$\\\hline

\ZZ{-6pt}{18pt }$(5).~x^2-y^2$ &$\frac{\partial}{\partial
x}\wedge\frac{\partial}{\partial y}+\pi_f$\\\hline

\ZZ{-6pt}{18pt}$(6).~x^2$ & $\frac{\partial}{\partial
x}\wedge\frac{\partial}{\partial y}+\pi_f$\\\hline
\end{tabular}
\begin{tabular}{|l|l|}\hline
\ZZ{-6pt}{18pt}quadratic function $f$& affine Poisson structure
\\\hline

\ZZ{-6pt}{18pt}$(7).~0$ & $\pi_{k,f}$\\\hline

\ZZ{-6pt}{18pt}$(8).~a(x^2+y^2)$& $\pi_{k,f}$  \\\hline

\ZZ{-6pt}{18pt}$(9)_1.~a(x^2-y^2),~a\neq\frac{1}{4}$& $\pi_{k,f}$
\\\hline

\ZZ{-6pt}{18pt}$(9)_2.~\frac{1}{4}(x^2-y^2),$ &
$\frac{\partial}{\partial z}\wedge\frac{\partial}{\partial
x}+\pi_{k,f}$\\\hline

\ZZ{-6pt}{18pt}$(10).~x^2$ & $\pi_{k,f}$\\\hline
\end{tabular}

 where $a> 0 $ is a constant.
\end{thm}

As an application of the classification of 4-dimensional Lie
algebras, we give an example to describe  conformal symplectic
structure of corresponding linear Jacobi structure obtained by the
decomposition of linear Poisson structures on $\mathbb R^4$. Pursue
this geometric approach of Jacobi structure is very interesting and
we have another paper to study it. For More details about Jacobi
structure, conformal symplectic structure and contact structure,
please see \cite{ki:local} and  \cite{li:varietes}.

\begin{ex}
{\rm Consider the {\rm4}-dimensional Lie algebra decided by $f=x^2$
and the derivation $
D=\left(\begin{array}{ccc}\frac{1}{2}&0&0\\0&\frac{1}{4}&0\\0&0&\frac{1}{4}\end{array}\right)$,
the corresponding linear Poisson structure is
$$
\pi=2x_1\frac{\partial}{\partial x_2}\wedge\frac{\partial}{\partial
x_3}+\frac{1}{2}x_1\frac{\partial}{\partial
x_1}\wedge\frac{\partial}{\partial
x_4}+\frac{1}{4}x_2\frac{\partial}{\partial
x_2}\wedge\frac{\partial}{\partial
x_4}+\frac{1}{4}x_3\frac{\partial}{\partial
x_3}\wedge\frac{\partial}{\partial x_4}
$$
and $k=D(\pi)=(0,0,0,1)$, so the corresponding Jacobi structure is
$$
(E,\Lambda)=(\frac{1}{3}\frac{\partial}{\partial
x_4},~2x_1\frac{\partial}{\partial
x_2}\wedge\frac{\partial}{\partial
x_3}+\frac{1}{6}x_1\frac{\partial}{\partial
x_1}\wedge\frac{\partial}{\partial
x_4}-\frac{1}{12}x_2\frac{\partial}{\partial
x_2}\wedge\frac{\partial}{\partial
x_4}-\frac{1}{12}x_3\frac{\partial}{\partial
x_3}\wedge\frac{\partial}{\partial x_4}).
$$
After some straightforward computations we have, if $x_1\neq 0$, the
character distribution is 4-dimensional and the 2-form $\Omega$,
inverse to the bi-vector $\Lambda$, is
$$
\Omega=\frac{1}{4}\frac{x_3}{x_1^2}dx_1\wedge
dx_2+\frac{1}{4}\frac{x_2}{x_1^2}dx_3\wedge
dx_1+\frac{6}{x_1}dx_4\wedge dx_1+\frac{1}{2x_1}dx_3\wedge dx_2.
$$
and $ \omega=i_E\Omega=\frac{2}{x_1}dx_1.$  So we have
$$
d\Omega=-\omega\wedge\Omega=\frac{1}{x_1^2}dx_1\wedge dx_2\wedge
dx_3.
$$
This shows that $(\Omega,\omega)$ is the the conformal symplectic
structure on the leaf.

If $x_1=0$, it is evident that the character distribution through
$(0,x_2^0,x_3^0,x_4^0)$ is 2-dimensional if  $x_2^0,x_3^0$ are not
zero at the same time and the leaf is just the half-plane decided by
$x_4$-axis and the point $(0,x_2^0,x_3^0,x_4^0)$ with $x_4$-axis
omitted; The character distribution through $(0,0,0,x_4^0)$ is
1-dimensional and the leaf is just $x_4$-axis.}  \qed
\end{ex}

\end{document}